\documentclass[reqno,fleqn,11pt]{amsart}
 \normalbaselines
\usepackage[frenchb,english]{babel}
\usepackage{amsmath}
\usepackage{graphicx}
\usepackage{amsfonts}
\usepackage{amssymb}%
\def\cal#1{\mathcal{#1}}
\def\Iso{\mathop{\hbox{Iso}}}
\def\BW{\mathop{\hbox{BW}}}

\def\det{\mathop{\hbox{det}}}
\def\tr{\mathop{\hbox{tr}}}

\def\G{\mathop{\hbox{Grass}}}
\def\GL{\mathop{\hbox{GL}}}

\def\span{\mathop{\hbox{{span}}}}

\def \L {{\cal L}}

\newtheorem{teor}{Theorem}[section]

\newtheorem{prop}[teor]{Proposition}
\newtheorem{corol}[teor]{Corollary}
\newtheorem{lemma}[teor]{Lemma}

\newcommand{\fdim}{\hspace*{\fill}$\Box$}
\newcommand{\dimostr}{{\bf Proof: }}

\newcommand{\real}{\Bbb{R}}

\newcommand{\complex}{\Bbb{C}}

\newcommand{\K}{K\"{a}hler}

\begin{document}

\noindent \centerline {\large\bf  SYMPLECTIC DUALITY OF SYMMETRIC
SPACES }



\vspace{0.5cm}

\centerline{\small Antonio J. Di Scala} \centerline{\small
Dipartimento di Matematica-Politecnico di Torino}
\centerline{\small Corso Duca degli Abruzzi 24, 10129
Torino-Italy}

\centerline{\small e-mail address: antonio.discala@polito.it}

\vspace{0.3cm} \centerline{\small and} \vspace{0.3cm}

\centerline{\small Andrea Loi} \centerline{\small Dipartimento di
Matematica e Informatica -- Universit\`{a} di Cagliari}
\centerline{\small Via Ospedale 72, 09124 Cagliari-- Italy}
\centerline{\small e-mail address: loi@unica.it}

\vskip 0.5cm

\centerline{\bf Abstract} {\small Let $(M, 0)  \subset
{\complex}^n$ be a complex $n$-dimensional Hermitian symmetric
space endowed  with the hyperbolic form $\omega_{B}$. Denote by
$(M^*, \omega_B^*)$ the compact dual of  $(M, \omega_B)$, where
$\omega_B^*$ is the Fubini--Study form on $M^*$. Our first result
is Theorem \ref{mainteor} where, with the aid of the theory of
Jordan triple systems, we construct an explicit diffeomorphism
$\Psi_M: M\rightarrow {\real}^{2n}={\complex}^n\subset M^*$
satisfying $\Psi_M^*(\omega_0)=\omega_B$ and
$\Psi_M^*(\omega_B^*)=\omega_0$. Amongst other properties of the
map $\Psi_M$, we also show that it takes (complete) complex and
totally geodesic submanifolds of $M$ through the origin to complex
linear subspaces of ${\complex}^n$. As a byproduct of the proof of
Theorem \ref{mainteor} we get an interesting characterization of
the Bergman form on a Hermitian space in terms of its restriction
to classical Hermitian symmetric spaces of noncompact type (see
Theorem \ref{secondteor} below).

\vskip 0.3cm


\noindent {\it{Subj.Class}}: 53D05, 58F06.}

\section{Introduction}
Dusa McDuff \cite{mc} proved a global version of  Darboux theorem
for a $n$-dimensional \K\ manifold with a pole $p$ such that its
radial curvature  is nonpositive, by showing that there exists a
diffeomorphism $\psi_p :M\rightarrow {\real}^{2n}={\complex}^n$
(depending on $p$) such that $\psi_p(p)=0$
$\psi_p^*(\omega_{0})=\omega$, where
$\omega_0=\sum_{j=1}^ndx_j\wedge dy_j$ is the standard symplectic
form on ${\real}^{2n}$. The interest for these kind of questions
comes, for example, after Gromov's discovery \cite{gr} of the
existence of exotic symplectic structures on ${\real}^{2n}$.
Eleonora  Ciriza \cite{cr2} (see also \cite{cr1} and \cite{cr3})
proves that the image $\psi_p(T)$ of any (complete)  complex and
totally geodesic  submanifold $T$ of  $M$ which contains the pole
is a complex linear subspace of ${\complex}^n$. It is worth
pointing out that McDuff's proof deals only with the existence
problem and the expression of the symplectomorphism $\psi_p$ is,
in general, very hard to find (see \cite{alcu} for the
construction of explicit symplectic coordinates on some complex
domains in ${\complex}^n$).

\vskip 0.3cm

\noindent In this paper we deal with the symplectic geometry of
Hermitian symmetric spaces of noncompact type. We are going to
regard such spaces as bounded symmetric domains $(M, 0)  \subset
{\mathcal M}$ centered in the origin
 of their associated Hermitian positive
Jordan triple system ${\mathcal M}$. Furthermore $M$ will be
equipped with the hyperbolic form. Let $(M^*, \omega_B^*)$ be the
compact dual of $(M, \omega_B)$ equipped with the Fubini--Study
form $\omega_B^*$ (see the next section for the definition of
hyperbolic and Fubini--Study form). We denote with the same symbol
the \K\ form $\omega_B^*$ on ${\mathcal M}$ obtained by the
restriction of  $\omega_B^*$ via the Borel embedding ${\mathcal M}
\subset M^*$. Finally, we denote by $HSSNT$ the space of all
Hermitian symmetric spaces of noncompact type $(M, 0)$ and by
${\mathcal P}$ the set of all diffeomorphisms $\psi: M\rightarrow
{\mathcal M}$, $M\in HSSNT$, such that $\psi (0)=0$.

\vskip 0.3cm

\noindent Our main result is the following theorem which
establishes a bridge among the symplectic geometry of HSSNT, their
duals and the theory of Jordan triple systems.

\begin{teor}\label{mainteor}
Under the above assumptions, the map
\begin{equation}\label{PsiM}
\Psi_M : M \rightarrow  {\mathcal M},\  z\mapsto B(z,
z)^{-\frac{1}{4}}z ,\end{equation}

has the following properties:

\begin{itemize}

\item[(D)] $\Psi_M$ is a  (real analytic) {\em diffeomorphism}
and its inverse $\Psi^{-1}_M$ is given by:
\begin{equation*}
\Psi^{-1}_M :  {\mathcal M} \rightarrow M,\ z\mapsto B(z, -
z)^{-\frac{1}{4}}z  \hspace{.1cm} ;
\end{equation*}

\item[(H)] The map $\Psi:HSSNT\rightarrow {\mathcal P}$ which
takes an $M\in HSSNT$  into the diffeomorphism $\Psi_M$ is  {\em
hereditary} in the following sense: for any $(T,0)
\stackrel{i}{\hookrightarrow} (M,0) $ complex and totally geodesic
embedded submanifold $(T,0)$  through the origin $0$, i.e.
$i(0)=0$ one has:
\begin{equation*}
\Psi_M {_ {|_T}} =  \Psi_T.
\end{equation*}
Moreover
\begin{equation}
\Psi_M (T)={\mathcal T}\subset {\mathcal M},
\end{equation}
where ${\mathcal T}$ is the Hermitian positive Jordan triple
system associated to $T$;

\item[(I)] $\Psi_M$ is a (non-linear) {\em interwining}  map w.r.t
the action of the isotropy group $K \subset \Iso(M)$ at the
origin, where $\Iso(M)$ is the group of isometries of $M$, i.e.
for every $\tau \in K$
\begin{equation*}
\Psi_M \circ \tau = \tau \circ \Psi_M ;
\end{equation*}

\item[(S)] $\Psi_M$ is a {\em symplectic duality}, i.e. the
following holds
\begin{equation}\label{s1}
\Psi_M^*(\omega_B^*)=\omega_0;
\end{equation}
\begin{equation}\label{s2}
\Psi_M^*(\omega_0)=\omega_B,
\end{equation}
where $\omega_0$ is the flat \K\ form on ${\mathcal M}$ (see
formula (\ref{flatform}) below for its definition).

\end{itemize}

\end{teor}

\vskip 0.3cm

\noindent From the point of view of inducing geometric structures
as in Gromov's programme \cite{gr} the importance of property (S)
relies on the existence of a smooth map (i.e. $\Psi_M$)  which is
a simultaneous symplectomorphism with respect to different
symplectic structures, i.e. satisfying (\ref{s1}) and (\ref{s2}).
(We refer the reader to \cite{c1nash} and \cite{connmet} and the
reference therein for the case of induction of different pairs
like symplectic forms and Riemannian metrics or connections and
Riemannian metrics). Notice also that property  (H) is exactly the
above mentioned property observed by Ciriza for the McDuff map,
namely the image of a totally geodesic submanifold (through the
origin) via the map $\Psi_M$ is sent to a complex linear subspace
of ${\mathcal M}$.

\vskip 0.3cm

\noindent The map $\Psi_M: M\rightarrow {\mathcal M}$ above was
defined, independently from the authors, by Guy Roos in
\cite{roos2} (see Definition VII.4.1 at p. 533). There (see
Theorem VII.4.3) he proved the analogous of (S) for volumes,
namely $\Psi_M^*(\omega_0^n)=\omega_B^n$ and
$\Psi_M^*((\omega_B^*)^n)=\omega_0^n$ ($n$ is the complex
dimension of $M$) which is, of course a corollary of (S).

The case where $M$ is the first Cartan domain $D_I[n]$ (namely the
dual of $\G _n ({\complex}^{2n})$) the map $\Psi_{M}$ was already
considered by John Rawnsley in a unpublished work of 1989, where
he proved property (\ref{s1}) for this case (see Section
\ref{basicexample} below for details). Actually the proof of (S)
for classical HSSNT (i.e. those Hermitian spaces which do not
contain exceptional factors in theirs De--Rham decomposition)
follows  from the property (S) for $D_I[n]$ (see  \ref{subclass}
below). Regarding the proof of (S) in the general case we present
here two proofs. The first one, presented in Section \ref{secgen},
is actually a \lq\lq partial proof'' since it is obtained by
assuming that one already knows that the symplectic form
$\Psi_M^*(\omega_0)$ and $(\Psi_M^{-1})^*(\omega_B^*)$ are of type
$(1,1)$. The second (and complete) proof, more algebraic in
nature, is due to Guy Roos. His proof  is, as it often occurs,
more or less an adaptation of the proof for matrices into the
language of Jordan triples and their operators.

\noindent There are two reasons of having included our (partial)
proof in this paper. First, it is of geometric nature and second
because, the techniques  employed, heuristically, have suggested
us how to attack and prove Theorem \ref{secondteor} below which
gives an interesting (and to the authors' knowledge unknown)
characterization of the Bergman metric in terms of its restriction
to classical HSSNT.

\vskip 0.3cm

\noindent The paper is organized as follows. In the next section
we collect some basic material about Hermitian positive Jordan
triple systems, Hermitian symmetric spaces and their dual. Section
\ref{basicexample} is dedicated to the proof of (D) and (S) of
Theorem \ref{mainteor} for the first Cartan domain. The result of
these sections are used in Section \ref{secclass} to prove (H) and
(I) of Theorem \ref{mainteor} in the general case and Theorem
\ref{mainteor} in the classical case. In Section \ref{secgen},
after recalling some basic facts on Jordan algebras we prove (D)
and (S) of Theorem \ref{mainteor} by reduction to the classical
case (property (S) is proved only in the hypothesis mentioned
above). Moreover at the end of this section  we state Theorem
\ref{secondteor} whose proof can be easily obtained by the same
method used in the proof of (S). Finally Section \ref{secroos}
contains Roos's proof of  (S). The paper ends with an appendix
containing two technical results on Hermitian positive Jordan
triple systems.

\vskip 0.3cm

\noindent We would like to thank Guy Roos for giving us the
opportunity of including his proof in this paper and for useful
discussions about his work on Jordan triple systems and for his
interest in ours.

\section{Jordan triple systems, Hermitian spaces of noncompact type and their compact dual}

\noindent
\subsection{Jordan triple systems and Hermitian spaces}

\noindent We briefly recall some standard material about Hermitian
symmetric spaces of noncompact type and Hermitian positive Jordan
triple systems. We refer to  \cite{roos2} for details, notations
and further results.

\vskip 0.3cm

\noindent An Hermitian Jordan triple system is a  pair $({\mathcal
M}, \{ ,\  ,\})$, where ${\mathcal M}$ is a complex vector space
and $\{ ,\  ,\}$ is a ${\real}$-trilinear map
\begin{equation*}
\{ ,\  ,\}:{\mathcal M}\times {\mathcal M}\times {\mathcal M}
\rightarrow {\mathcal M}, (u, v, w)\mapsto \{u, v, w\}
\end{equation*}
which is ${\complex}$-bilinear and symmetric in $u$ and $w$,
${\complex}$-antilinear in $v$ and such that the following {\em
Jordan identity} holds:
\begin{equation}\label{jordan}
\{x, y, \{u, v, w\}\}-\{u, v, \{x, y, w\}\}= \{\{x, y, u\}, v,
w\}-\{u, \{v, x, y\}, w\}.
\end{equation}

For $u, v\in {\mathcal M}$, denote by $D(u, v)$ the operator of
${\mathcal M}$
\begin{equation*}
D(u, v)(w)=\{u, v, w\}.
\end{equation*}
An Hermitian Jordan triple system is called  {\em positive} if the
Hermitian form $(u, v)\mapsto \tr D(u, v)$ is positive definite.
In the sequel we will write $HPJTS$ to denote an Hermitian
positive Jordan triple system. We also denote (with a slight abuse
of notation) by $HPJTS$ the set of all Hermitian positive Jordan
triple systems on a fixed complex vector space ${\mathcal M}$. An
$HPJTS$ is always {\em semi-simple}, that is a finite family of
simple subsystems with component-wise triple product. A $HPJTS$ is
called {\em simple} if it is not the product of two non-trivial
Hermitian positive Jordan triple subsystems. The {\em quadratic
representation} $Q:{\mathcal M}\rightarrow End_{\real}({\mathcal
M})$ is  defined by
\begin{equation*}
2Q(u)(v)=\{u, v, u\}, \ u, v\in {\mathcal M}.
\end{equation*}
The {\em Bergman operator} is defined by
\begin{equation}\label{bergmanoperator}
B(u, v)=\operatorname{id}-D(u, v)+Q(u)Q(v),
\end{equation}
where $\operatorname{id}: {\mathcal M} \rightarrow {\mathcal M}$
denotes the identity map of ${\mathcal M}$.

\vskip 0.3cm

\noindent An element $c\in {\mathcal M}$ is called {\em tripotent}
if $\{c,c,c\}=2c$. Two tripotents $c_1$ and $c_2$ are called {\em
orthogonal} if $D(c_1, c_2)=0$. A non zero tripotent $c$ is called
{\em primitive} if it is not the sum of non-zero orthogonal
tripotents. Due to the positivity of the Jordan triple system
${\mathcal M}$, each element $x\in {\mathcal M}$
has a unique \emph{spectral decomposition }%
\begin{equation}
x=\lambda_{1}c_{1}+\lambda_{2}c_{2}+\cdots+\lambda_{p}c_{p}, \label{eqr1}%
\end{equation}
where $\lambda_{1}>\lambda_{2}>\cdots>\lambda_{p}>0$ and $\left(  c_{1}%
,\ldots,c_{p}\right)  $ is a system of mutually orthogonal
tripotents. Moreover, each $x\in {\mathcal M}$ may also be written
\begin{equation}
x=\lambda_{1}c_{1}+\lambda_{2}c_{2}+\cdots+\lambda_{r}c_{r}, \label{eqr2}%
\end{equation}
with $\lambda_{1}\geq\lambda_{2}\geq\cdots\geq\lambda_{r}\geq0$
and $\left( c_{1},\ldots,c_{r}\right)  $ is a \emph{frame} (that
is a maximal system of mutually orthogonal tripotents). The
decomposition (\ref{eq2}) is also called {\em spectral
decomposition}; it is unique only for elements $x$ of maximal rank
$r$, which form a Zariski dense open subset of ${\mathcal M}$.

There exist polynomials $m_1,\dots ,m_r$ on ${\mathcal M}\times
\overline{\mathcal M}$, homogeneous of respective bidegrees $(1,
1),\dots ,(r, r)$, such that for $x\in {\mathcal M}$, the
polynomial
\begin{equation*}
m(T, x, y)=T^r-m_1(x, y)T^{r-1}+\cdots +(-1)^rm_r(x, y)
\end{equation*}
satisfies
\begin{equation*}
m(T, x, x)=\prod_{i=1}^{r}(T-\lambda_i^2),
\end{equation*}
where $x$ is the spectral decomposition of $x=\sum \lambda_jc_j$.

\noindent The inohomogeneous polynomial
\begin{equation*}
N(x, y)=m(1, x, y)
\end{equation*}
is called the generic norm.

 Denote by $\mathcal{N}$ and $\mathcal{N}_{*}$ the
associated functions%
\begin{equation}\label{N}
\mathcal{N}(z)=N(z,z)=1-m_{1}(x,x)+\cdots+(-1)^{k}m_{k}(x,x)+\cdots
+(-1)^{r}m_{r}(x,x),
\end{equation}
and
\begin{equation}\label{N*}
\mathcal{N}_{*}(z)=N(z,-z)=1+m_{1}(x,x)+\cdots+m_{k}(x,x)+\cdots+m_{r}(x,x).
\end{equation}

\subsection{HSSNT associated to HPJTS}
M. Koecher (\cite{ko1}, \cite{ko2}) discovered that to every HPJTS
$({\mathcal M}, \{ ,\  ,\})$ one can associate an Hermitian
symmetric space of noncompact type, i.e. a bounded symmetric
domain $(M, 0)$ centered at the origin $0\in {\mathcal M}$. The
domain $(M, 0)$ is defined as the connected component containig
the origin of   the set of all $u\in {\mathcal M}$ such that $B(u,
u)$ is positive definite with respect to the Hermitian form $(u,
v)\mapsto \tr D(u, v)$. The Bergman form $\omega_{Berg}$ on $M$ is
given by
\begin{equation*}
\omega_{Berg} =-\frac{i}{2\pi}\partial\bar\partial\log\det B.
\end{equation*}
The {\em hyperbolic metric} $\omega_B$ (which appears in the
statement of our Theorem \ref{mainteor}) is given by:
\begin{equation}\label{bergmanform}
\omega_B = -\frac{i}{2\pi}\partial\bar\partial\log{\mathcal N},
\end{equation}
where ${\mathcal N}(z)$ is given by (\ref{N}). If $M$ is
irreducible, or equivalently ${\mathcal M}$ is simple, then  $\det
B={\mathcal N}^{g}$, where $g$ is the genus of $M$, and hence, in
this case, $\omega_B =\frac{\omega_{Berg}}{g}$. Observe that in
the rank one case, that is when $M$ is the complex Hermitian ball,
the form $\omega_B$ is the standard hyperbolic metric (cfr.
formula (\ref{hypgras}) in the next section ).

The HPJTS $({\mathcal M}, \{ ,\  ,\})$ can be recovered by its
associated Hermitian symmetric space of noncompact type $(M, 0)$
by defining ${\mathcal M}=T_0M$ (the tangent space to the origin
of $M$) and
\begin{equation}\label{trcurv}
\{u, v, w\}=-\frac{1}{2}\left(R_0(u, v)w+J_0R_0(u, J_0v)w\right),
\end{equation}
where $R_0$ (resp. $J_0$) is the curvature tensor of the Bergman
metric (resp. the complex structure) of $M$ evaluated at the
origin. For more informations on the correspondence between
$HPJTS$ and $HSSNT$ we refer to p. 85 in Satake's book \cite{sa}
(see also \cite{loos}). We refer also  to \cite{Ber} for some deep
implications of formula (\ref{trcurv}).

\subsection{Totally geodesic submanifolds of HSSNT}
In the proof of our theorems  we  need the following result whose
proof follows by equality (\ref{trcurv}) and the well-known
correspondence between totally geodesic submanifolds and Lie
triple systems (see Theorem 4.3 p. 237 in \cite{kn}).
\begin{prop} \label{geo}
Let $(M,0)$ be a HSSNT with origin $0 \in M$ and let ${\mathcalÊ
M}$ be its associated $HPJTS$. Then there exists a one to one
correspondence between (complete) complex totally geodesic
submanifolds and sub-$HPJTS$ of ${\mathcal M}$. This
correspondence sends $(T,0) \subset (M,0)$ to  ${\mathcal
T}\subset {\mathcal M}$, where ${\mathcal T}$ denotes the HPJTS
associated to $T$.
\end{prop}

\subsection{The compact dual of an  $HSSNT$.}
Let $(M^*, \omega_B^*)$ be  the compact dual of $(M, \omega_B)$.
It is well known (see e.g. \cite{ta}) that $(M^*, \omega_B^*) $ is
a compact homogeneous simply-connected \K\ manifold  and $M^*$
admits a holomorphic embedding $\BW:M^*\rightarrow {\complex}P^N$
(the {\em Borel--Weil embedding}) into a $N$-dimensional complex
projective space satisfying $\BW^*(\omega_{FS})=\omega_B^*$, where
$\omega_{FS}$ denotes the Fubini--Study  form on
${\complex}P^{N}$, namely the \K\ form which,  in the homogeneous
coordinates $[z_0,\dots z_N]$ on ${\complex}P^{N}$, is given by
$\omega_{FS}=\frac{i}{2\pi}
\partial\bar\partial\log (|z_0|^2+\cdots +|z_N|^2)$.
In our Theorem \ref{mainteor} and in the sequel we will call
$\omega_B^*$ the {\em Fubini--Study form} on $M^*$. In order to
write its local expression, let $p\in M^*$ and assume, without
loss of generality, that $\BW(p)=[1,0,\dots , 0]\in
{\complex}P^N$. Let $H_p\subset {\complex}P^N$ be the hyperplane
at infinity corresponding to the point $BW(p)$ and set
$Y_p=\BW^{-1}(H_p)$. One can prove (see  \cite{wolf}) that
$M^*\setminus Y_p$ is biholomorphic to ${\mathcal M}=T_0M$ (the
$HPJTS$ associated to $M$). Moreover, under the previous
biholomorphism, $p$ can be made to correspond to the origin $0\in
M$. Hence we have the following inclusions $M \subset {\mathcal M}
\subset M^*$ and one can prove that  the restriction to ${\mathcal
M}$ of the \K\ form $\omega_{B}^*$ is given by:
\begin{equation}\label{sympdual}
\omega_B^* =\frac{i}{2\pi}\partial\bar\partial\log\det {\mathcal
N}_*,
\end{equation}
where ${\mathcal N}_*(z)$ is given by (\ref{N*}) (see also
\cite{diastherm} for the relations between the two \K\ forms
$\omega_B$ and $\omega_B^*$).

\subsection{The flat form on ${\mathcal M}$}
The {\em flat \K\ form} on ${\mathcal M}$ is defined by
\begin{equation}\label{flatform}
\omega_0=\frac{i}{2\pi}\partial\bar\partial m_1(x, x),
\end{equation}
where $m_1(x, x)$ is the polynomial appearing in (\ref{N}).
Observe that if ${\mathcal M}$ is simple then $\tr D(x, y)=gm_1(x,
y)$ and hence $\omega_0=\frac{i}{2g\pi}\partial\bar\partial D(x,
x)$. Notice also that in the rank-one case $\omega_0$ is (up to
the factor $2\pi$) the standard Euclidean form on ${\mathcal
M}={\complex}^n$ (cfr. formula (\ref{flatgras}) below).

\section{The proof of (D) and (S) of Theorem \ref{mainteor} for the first Cartan's domain}\label{basicexample}
Let $M=D_1[n]$ be the complex noncompact dual of
$M^*=G_n({\complex}^{2n})$, where $G_n({\complex}^{2n})$ is the
complex Grassmannian of complex $n$ subspaces of
${\complex}^{2n}$. In its realization as a bounded domain,
$D_1[n]$ is given by
\begin{equation}
D_1[n]=\{Z\in M_{n}({\complex})|\ I_n-ZZ^*>>0\},
\end{equation}

The triple product on ${\complex}^{n^2}$ making it an  HPJTS is
given by
\begin{equation}\label{tripgrass}
\{U, V, W\}=UV^*W+WV^*U, \ U, V, W\in M_n({\complex}).
\end{equation}
Hence the Bergman operator is given by
\begin{equation*}
B(U, V)W=(I_n-UV^*)W(I_n-V^*U)
\end{equation*}
A simple computation shows that the hyperbolic form and the
duality   map $\Psi_M:D_I[n]\rightarrow
M_n({\complex})={\complex}^{n^2}$ are given  by:
\begin{equation}\label{hypgras}
\omega_{B}= -\frac{i}{2\pi}\partial\bar\partial\log \det
(I_n-ZZ^*).
\end{equation}
\begin{equation}\label{rawnmap}
\Psi_M(Z)=(I_n-ZZ^*)^{-\frac{1}{2}}Z,
\end{equation}
respectively. In this case the Borel--Weil embedding is given by
the Pl\"{u}cker embedding $G_n({\complex}^{2n})\hookrightarrow
{\complex}P^N, N=(\begin{array}{cc}
2n\\
n
\end{array})-1$
and the local expression (\ref{sympdual}) of $\omega_B^*$ on
${\cal M}={\complex}^{2n}$ reads as
\begin{equation}\label{grloc}
\omega_{B}^*= \frac{i}{2\pi}\partial \bar\partial\log\det (I_n +
XX^*),
\end{equation}
with $X\in M_n(\complex)$. Moreover the flat \K\ form
(\ref{flatform}) is given by
\begin{equation}\label{flatgras}
\omega_0=\frac{i}{2\pi}\partial\bar\partial\log\tr (ZZ^*).
\end{equation}
By using the equality
\begin{equation*}
XX^*(I_n+ XX^*)^{\frac{1}{2}}=(I_n+ XX^*)^{\frac{1}{2}}XX^*
\end{equation*}
it is easy to verify  that the map
\begin{equation}
\Phi_M: {\complex}^{n^2}\rightarrow D_1[n], X\mapsto
(I_n+XX^*)^{-\frac{1}{2}}X
\end{equation}
is the inverse of $\Psi_M$.

We are now ready to prove (S), namely the equalities
\begin{equation}\label{eq1}
\Psi_M^*(\omega_0)=\omega_B
\end{equation}
\begin{equation}\label{eq2}
\Psi_M^*(\omega_B^*)=\omega_0.
\end{equation}
As we already pointed out, the proof of the  equation (\ref{eq1}),
is due to J. Rawnsley (unpublished). Here we present his proof.
First of all observe that we can write
\begin{eqnarray*}
{\omega_B}&=&-\frac{i}{2\pi}\partial\bar\partial\log\det
(I_n-ZZ^*)=
\frac{i}{2\pi}\operatorname*{d}\partial\log\det (I_n-ZZ^*)\\
&=&\frac{i}{2\pi}\operatorname*{d}\partial\tr\log (I_n-ZZ^*)=\frac{i}{2\pi}\operatorname*{d}\tr\partial\log (I_n-ZZ^*)\\
&=&-\frac{i}{2\i}\operatorname*{d}\tr[Z^*(I_n-ZZ^*)^{-1}dZ],
\end{eqnarray*}
where we use the decomposition $\operatorname*{d}=\partial
+\bar\partial$ and the identity $\log\det A=\tr\log A$. By
substituting $X=(I_n-ZZ^*)^{-\frac{1}{2}}Z$ in the last expression
one   gets:
\begin{equation*}
-\frac{i}{2\pi}\operatorname*{d}\tr[Z^*(I_n-ZZ^*)^{-1}dZ] =
-\frac{i}{2\pi}\operatorname*{d}\tr(X^*dX)+\frac{i}{2\pi}
\operatorname*{d}\tr\{X^*d[(I_n-ZZ^*)^{-\frac{1}{2}}]Z\}.
\end{equation*}

Observe now that
$-\frac{i}{2\pi}\operatorname*{d}\tr(X^*dX)=\omega_0$ and the
$1$-form $\tr[X^*\operatorname*{d}(I_n-ZZ^*)^{-\frac{1}{2}}Z]$ on
${\complex}^{n^2}$ is exact being equal to $\operatorname*{d}\tr
(\frac{C^2}{2}-\log C)$, where $C=(I_n-ZZ^*)^{-\frac{1}{2}}$.
Therefore $\omega_B$ in the $X$-coordinates equals $\omega_0$ and
this concludes the proof of equality (\ref{eq1}).

The proof of  (\ref{eq2}) follows the same line. Indeed, by
(\ref{grloc}) we get
\begin{eqnarray*}
\omega_{B}^*&=&\frac{i}{2\pi}\partial\bar\partial\log\det
(I_n+XX^*)=
-\frac{i}{2\pi}\operatorname*{d}\tr\partial\log (I_n+XX^*)\\
&=&-\frac{i}{2\pi}\operatorname*{d}\tr[X^*(I_n+XX^*)^{-1}dX]
\end{eqnarray*}

By substituting $Z=(I_n+XX^*)^{-\frac{1}{2}}X$ in the last
expression one   gets:
\begin{eqnarray*}
-\frac{i}{2\pi}\operatorname*{d}\tr[X^*(I_n-XX^*)^{-1}dX] &=&
-\frac{i}{2\pi}\operatorname*{d}\tr(Z^*dZ)
+\frac{i}{2\pi}\operatorname*{d}\tr\{Z^*d[(I_n+XX^*)^{-\frac{1}{2}}]X\}\\
&=&\omega_0+\frac{i}{2\pi}{\operatorname*{d}}^2\tr(\log
D-\tr\frac{D^2}{2})=\omega_0,
\end{eqnarray*}
where $D=(I_n+XX^*)^{-\frac{1}{2}}$ and this concludes the proof
of (S) for $D_I[n]$.

\section{Proof of (H) and (I) and the proof of Theorem \ref{mainteor}
for classical domains}\label{secclass}

Let $(M, 0)$ be any HSSNT. Since the map $\Psi_M$ depends  only on
the triple product $\{,\ ,\}$ properties (H) is a straightforward
consequence of Proposition \ref{geo} above. Let ${\mathcal M}$ be
the HPJTS associated to $M$. As usual, let us write $M=G/K$, where
$G$ is the isometry group of $M$ and $K\subset G$ is the (compact)
isotropy subgroup of the origin $0\in M$. Due to a theorem of E.
Cartan (see \cite{mok}, p. 63) the group $K$ consists entirely of
linear transformations, i.e. $K\subset\GL ({\mathcal M})$. In
order to prove (I) of Theorem \ref{mainteor}, observe that the
Bergman operator associated to ${\mathcal M}$ is invariant by the
group of isometry of $M$, namely for every isometry $\tau\in K$
\begin{equation*}
B(\tau (u), \tau (v))(\tau (w))=\tau\left(B(u, v)(w)\right),
\forall u, v, w\in {\mathcal M},
\end{equation*}
which implies that
\begin{equation*}
B(\tau(z) \, , \, \tau(z))^{-1/4} (\cdot) =
\tau(B(z,z)^{-1/4}(\tau^{-1}(\cdot))), \ \forall z\in M.
\end{equation*}
Hence
\begin{equation*}
\Psi_M\circ\tau =\tau\circ \Psi_M
\end{equation*}
for all $\tau\in K$ and we are done. \fdim

\vskip 0.3cm

\noindent
\subsection{Proof of Theorem \ref{mainteor} for classical HSSNT}\label{subclass}
Observe that since now we  have proved properties (H), (I)  for
any HSSNT and property (D) and (S)  for $D_I[n]$. Let $(M, 0)$ be
a classical HSSNT and let ${\mathcal M}$ be its associated
$HPJTS$. It is well-known that $(M, 0)$ can be complex and totally
geodesic embedded  into $D_I[n]$, for $n$ sufficiently large. We
can assume that this embedding takes the origin $0\in M$ to the
origin $0\in D_I[n]$. Therefore, by Proposition \ref{geo}, the
HPJTS ${\mathcal M}$ is a sub-HPJTS of $({\complex}^{n^2},
\{,,\})$. Hence properties (D), (S)  for $M$ are consequences of
property (H) and the fact (proved in Section \ref{basicexample})
that these properties hold true for $D_I[n]$.

\section{Proof of Theorem \ref{mainteor} by reduction to the classical case}
\label{secgen} Let $(M, 0)$ be an exceptional  $HSSNT$. Properties
(H) and (I) were proved in the previous section. In this section
we prove (D) and (S) In order to prove (D) and (S), we pause to
prove Lemmata \ref{extension}, \ref{inseticida} and its Corollary
\ref{corolinseticida} below which will be the bridges between the
exceptional case and the classical one. Before  stating and
proving these lemmata let us briefly recall the concept of {\em
Jordan algebras} (see e.g. \cite{roos2} for details). A Jordan
algebra (over $\real$ or ${\complex}$) is a (real or complex)
vector space ${\cal A}$ endowed with a commutative bilinear
product
\begin{equation*}
\circ :{\cal A}\times {\cal A}\rightarrow {\cal A}, (a, b)\mapsto
a\circ b
\end{equation*}
satisfying the following identity:
\begin{equation*}
a\circ(a^2\circ b)=a^2\circ (a\circ b), \forall a, b\in {\cal A},
\end{equation*}
where $a^2=a\circ a$. Given a Jordan algebra ${\cal A}$ over
${\complex}$ the triple product given by
\begin{equation*}
\{x, y, z\}=2\left((x\circ \bar y)\circ z+(z\circ \bar y)\circ x -
(x\circ z)\circ \bar y\right )
\end{equation*}
defines a structure of Jordan triple system  on ${\cal A}$ (cfr.
Proposition II.3.1 at page 459 and formula (6.18) at page 514 in
\cite{roos2}). Not all HPJTS ${\mathcal M}$ arises form a Jordan
algebra. If this happens the HSSNT associated to ${\mathcal M}$ is
called of {\em tube type}. Nevertheless we have the following
result

\begin{lemma} \label{extension}
Let $(M,0)$ be a HSSNT and let ${\cal M}$ be its associated HPJTS.
Then, there exists a HSSNT $(\widetilde{M},0)$ such that:

\begin{itemize}

 \item[(i)] $(M,0) \hookrightarrow (\widetilde{M},0)$ complex and
 totally geodesic embedded,

\item[(ii)] The HPJTS $\widetilde{\cal M}$ associated to
$(\widetilde{M},0)$ arises from a Jordan Algebra.

\end{itemize}

\end{lemma}
\dimostr Assume first that $M$ is irreducible. If $M$ is of
classical type take a suitable $n$ and  a complex and totally
geodesic embedding $(M, 0)\hookrightarrow D_I[n]$. The HPJTS
associated to $D_I[n])$ comes from a Jordan algebra and so, by
Proposition \ref{geo}, the lemma is proved for classical HSSNT. If
$M$ is of exceptional type, it is known (see Appendix in
\cite{roosall}) that the HPJTS  associated to $E_6$ (the
exceptional HSSNT of dimension $16$) is a sub-HPJTS of the HPJTS
$(H_3 ({\cal O}_{\complex}), \{,,\})$ of dimension $27$ associated
to the exceptional HSSNT $E_7$. Since $(H_3 ({\cal O}_{\complex}),
\{,,\})$ arises from a Jordan algebra (i.e. $E_7$ is of tube type,
see Appendix in \cite{roosall}), the prove of the lemma follows
again by Proposition \ref{geo}.

For a reducible HSSNT one simply takes the product of the Jordan
algebras associated to each factor. \fdim

\begin{lemma}\label{inseticida}
Let $M$ be a HSSNT. Let $p$ be a point of $M$,  $a, b\in
T_pM={\cal M}$ be two non-zero vectors and $\pi\subset T_pM$ be
the real plane generated by these vectors. Then there exists a
classical HSSNT $C \hookrightarrow M$ complex and totally
geodesically emdedded in $M$ passing trough $p$ such that $\pi
\subset T_pC$.
\end{lemma}
\dimostr Without loss of generality we can assume that  $p$ is
equal the origin $0$ of $M$. Consider the Jordan subalgebra ${\cal
C}_{ab}\subset \widetilde{\cal M}$ generated by $a$ and $b$, where
$\widetilde{\cal M}$ is the Jordan algebra given by  the previous
lemma. A deep result due to Jacobson-Shirsov \cite{JaPa} asserts
that this Jordan algebra is special, namely its associated HSSNT
is of classical type. Therefore, by (i) of the previous lemma, the
HSSNT $C\hookrightarrow M\hookrightarrow \widetilde M$ associated
to the HPJTS ${\cal C}_{ab}\cap {\cal M}\subset{\cal M}$ satisfies
the desired properties. \fdim

\begin{corol} \label{corolinseticida} Let $M$ be a HSSNT.
Let $p \in M$ be a point of $M$ and let $\gamma$ be a geodesic of
$M$ passing trough $p$. Then, for any complex line $\L \subset
T_pM$ there exists a classical HSSNT $C \hookrightarrow M$ complex
and totally geodesically emdedded in $M$ passing trough $p$ such
that $\gamma \subset C$ and $\L \subset T_pC$.
\end{corol}
\dimostr Let $a\in T_pM$ be a non zero vector  tangent to $\gamma$
and let $b\in T_pM$  any nonzero vector of ${\cal L}$. The desired
$C$ is then obtained by applying Lemma \ref{inseticida} to these
vectors. \fdim

\subsection{\bf Proof of (D)}
Let $p \in (M,0)$ be any point and let $v_p \in T_pM$ any vector
tangent to $M$. Then, by Lemma \ref{inseticida} we can construct a
classical (complete) complex totally geodesic $(C,0)
\hookrightarrow (M,0)$ such that $p \in N$ and $v_p \in T_pC$. By
property $(H)$ we know that $\Psi_M |_C$ is a diffeomorphism and
therefore we get that $(d \Psi_M)_p$ is bijective and therefore by
the inverse function theorem $\Psi_M$ is a local diffeomorphism.
In order to prove the injectivity of the $\psi_M$ let $p,q \in M$
and let $C \hookrightarrow M$ be  a classical HSSNT containing
$0,p,q \in C$, whose existence is guaranteed by Lemma
\ref{inseticida}. Then, property $(H)$ implies that $\Psi_M$ is
$1-1$.

In order to prove the surjectivity of $\Psi_M$, let $q \in {\cal
M}$ be an arbitrary point. We can assume that $q \neq 0$, since
$\Psi_M(0)=0$. We have to show that there exist $p \in M$ such
that $\Psi_M(p)=q$. Let now $\gamma$ be the $1$-dimensional real
subspace of ${\cal M}$ generated by $q$. Observe that regarding $M
\subset {\cal M}$ as a bounded domain, it follows that
$\widetilde{\gamma} = \gamma \cap M$ is (as a subset) a geodesic
of $M$. Let $C$ be a classical totally geodesic complex
submanifold of $M$ containing $\widetilde{\gamma}$ given by Lemma
\ref{inseticida}. Notice that by Proposition \ref{geo} the point
$q$ belongs to ${\cal C}$ (the HPJTS Êassociated to $C$). Since
$C$ is classical, we know (by Section \ref{secclass}) that  there
exists $p \in C \subset M$ such that $q = \Psi_C(p) = \Psi_M(p)$,
(last identity follows again by property (H)), and we are done.
Finally, note that the expression for the inverse of $\Psi_M$ also
follows from property $(H)$ by restriction to classical (complete)
complex totally geodesic submanifolds. \fdim

\subsection{\bf Proof of (S)
under the assumptions that $\Psi_M^*(\omega_0)$ and
$(\Psi_M^{-1})^*(\omega_B^*)$ are of type $(1, 1)$} We only give a
proof of (\ref{s1}) since (\ref{s2}) is obtained in a similar
manner by applying the following argument to the map
$\Psi_M^{-1}$.

First of all notice that if we set
$\omega_{\Psi_M}=\Psi_M^*(\omega_0)$ equality (\ref{s1}) is
equivalent to the validity of the following two equations
\begin{equation}\label{ns1}
(\omega_{\Psi_M})_p(u, Ju)=(\omega_B)_p(u, Ju),
\end{equation}
\begin{equation}\label{ns2}
(\omega_{\Psi_M})_p(Ju, Jv)=(\omega_B)_p(Ju, Jv),
\end{equation}
for all $p\in M$, $u, v\in T_pM$, where $J$ denotes the almost
complex structure of $M$ evaluated at the point $p$. The second
equation, namely (\ref{ns2}), is precisely our assumption that
$\Psi_M^*(\omega_0)$ is of type $(1, 1)$.

Thus it remains to prove (\ref{ns1}). Fix $p\in M$ and $u\in
T_pM$. Consider the complex line ${\cal
L}=\span_{\complex}(u)\subset T_pM$ and a classical complex and
totally geodesic submanifold $(C, 0)\hookrightarrow (M, 0)$ such
that ${\cal L}\subset T_pC$ (whose existence is guaranteed by
Corollary \ref{corolinseticida}). If we denote by $\omega_{B, C}$
and $\omega_{0,{\cal C}}$ the hyperbolic form on $C$ and the flat
\K\ form on ${\cal C}$ (the HPJTS associated to $C$) we get:
\begin{equation*}
(\omega_{\Psi_M})_p(u, Ju)=(\Psi_C^*(\omega_{0, {\cal C}}))_p (u,
Ju)=(\omega_{B, C})_p(u, Ju)=(\omega_B)_p(u, Ju),
\end{equation*}
where the first and third  equalities follow by the hereditary
property (H) and the fact that the embedding $(C,
0)\hookrightarrow (M, 0)$ is a \K\ embedding while the second
equality is true since $C$ is of classical type (and hence
$\Psi_C^*(\omega_{0, {\cal C}})=\omega_{B, C}$ by the results of
Section \ref{secclass}). \fdim

\vskip 0.3cm

\noindent As a byproduct of the previous proof one immediately
gets the following theorem.

\begin{teor}\label{secondteor}
Let $(M, 0)$ be a HSSNT equipped with its Bergman form
$\omega_{{Berg},M}$. Let  $\omega$ be a two form of type $(1, 1)$
on $M$. Assume that the restriction of $\omega$ to all classical
HSSNT $(C, 0)$ passing through the origin equals the Bergman form
of $C$. Then $\omega =\omega_{{Berg}, M}$.
\end{teor}

\section{Roos' proof of properties (D) and (S)}\label{secroos}

Let ${\mathcal M}$ be a HPJTS of rank $r$. In this section we
denote by $\Psi =\Psi_M: M\rightarrow {\mathcal M},z\mapsto
B(z,z)^{-1/4}z$ the duality map. Let
\begin{equation*}
z=\lambda_{1}c_{1}+\lambda_{2}c_{2}+\cdots+\lambda_{r}c_{r}
\end{equation*}
be a spectral decomposition of $z\in M$. As (see \cite{roos2},
Proposition V.4.2, (5.8))
\begin{eqnarray*}
B(z,z)c_{j}  &  =&\left(  1-\lambda_{j}^{2}\right)  ^{2}c_{j},\\
D(z,z)c_{j}  &  =&2\lambda_{j}^{2}c_{j},
\end{eqnarray*}
we have
\begin{equation}
\Psi(z)  =\sum_{j=1}^{r}\frac{\lambda_{j}}{\left(  1-\lambda_{j}%
^{2}\right)  ^{1/2}}c_{j}
\end{equation}
and
\begin{equation}
\left(  \operatorname{id}-\frac{1}{2}D(z,z)\right)  c_{j}
=\left( 1-\lambda_{j}^{2}\right)  c_{j}.
\end{equation}
Thus
\begin{equation}
\Psi(z)=\left(  \operatorname{id}-\frac{1}{2}D(z,z)\right)
^{-1/2}z=\left( \operatorname{id}-z\square z\right)  ^{-1/2}z
\label{eq11},
\end{equation}
where we use the operator
\begin{equation}
z\square z=\frac{1}{2}D(z,z). \label{eq15}%
\end{equation}

From the previous equation, it is easily seen that $\Psi$ is
bijective and that the inverse map $\Psi^{-1}:{\cal M}\rightarrow
M$ is given by
\begin{equation}
\Psi^{-1}(u)=\sum_{j=1}^{r}\frac{\mu_{j}}{\left(
1+\mu_{j}^{2}\right)
^{1/2}}c_{j}, \label{eq5}%
\end{equation}
if $u=\sum_{j=1}^{r}\mu_{j}c_{j}$ is the spectral decomposition of
$u\in V$.
The relation (\ref{eq5}) is equivalent to%
\begin{equation}
\Psi^{-1}(u)=B(u,-u)^{-1/4}u\qquad(u\in V), \label{eq6}%
\end{equation}
so that $\Psi$ is a diffeomorphism. Therefore (D) in Theorem
\ref{mainteor} is proved.

\vskip 0.3cm

\noindent In order to prove (S) of Theorem \ref{mainteor}, set
$p_{1}(z)=m_{1}(z,z).$ We then have
$\overline{\partial}p_{1}=m_{1}(z,\operatorname*{d}z)$ and

\begin{eqnarray}
\Psi^{\ast}\left(  \overline{\partial}p_{1}\right)   &
=&m_{1}\left(
\Psi(z),\operatorname*{d}\Psi(z)\right) \nonumber\\
&  =&m_{1}\left(  \Psi(z),\left(  \operatorname*{d}\left(  \operatorname{id}%
-z\square z\right)  ^{-1/2}\right)  z\right) \nonumber\\
&  +&m_{1}\left(  \Psi(z),\left(  \left(
\operatorname{id}-z\square z\right)
^{-1/2}\right)  \operatorname*{d}z\right),  \label{eq17b}%
\end{eqnarray}
where we have used the identity
\[
\operatorname*{d}\Psi(z)=\left(  \operatorname*{d}\left(  \operatorname{id}%
-z\square z\right)  ^{-1/2}\right)  z+\left(
\operatorname{id}-z\square z\right)  ^{-1/2}\operatorname*{d}z.
\]

As $z\square z$ is self-adjoint w.r. to the Hermitian metric
$m_{1}$, we have
\begin{eqnarray}
 m_{1}\left(  \Psi(z),\left(  \left(  \operatorname{id}-z\square z\right)
^{-1/2}\right)\operatorname{d}  z\right) \nonumber &
=&m_{1}\left(  \left(  \operatorname{id}-z\square z\right)
^{-1/2}z,\left( \left(  \operatorname{id}-z\square z\right)
^{-1/2}\right)  \operatorname*{d} z\right)\\ \nonumber
&  =&m_{1}\left(  \left(  \operatorname{id}-z\square z\right)  ^{-1}%
z,\operatorname*{d}z\right).  \label{eq17a}%
\end{eqnarray}
If $z=\lambda_{1}c_{1}+\lambda_{2}c_{2}+\cdots+\lambda_{r}c_{r}$
is a spectral decomposition of $z\in {\mathcal M}$, we have
\begin{equation}
\left(  \operatorname{id}-z\square z\right)  ^{-1}z=\sum_{j=1}^{r}%
\frac{\lambda_{j}}{1-\lambda_{j}^{2}}c_{j}=z^{z}, \label{eq14}%
\end{equation}
where $z^{z}$ denotes the \emph{quasi-inverse} in the Jordan
triple system ${\mathcal M}$.

Using (\ref{eq12}), (\ref{eq14}) and Lemma \ref{lemmaa1} in the
Appendix, we get the last term in (\ref{eq17b}), namely
\begin{equation}
m_{1}\left(  \Psi(z),\left(  \left(  \operatorname{id}-z\square
z\right) ^{-1/2}\right)  \operatorname*{d}z\right)
=-\frac{\overline{\partial
}\mathcal{N}}{\mathcal{N}}. \label{eq16}%
\end{equation}

Applying Lemma \ref{lemmaa2} in the Appendix, we get

\begin{eqnarray*}
m_{1}\left(  \Psi(z),\left(  \operatorname*{d}\left(  \operatorname{id}%
-z\square z\right)  ^{-1/2}\right)  z\right) &=&m_{1}\left(
\Psi(z),\frac{1}{2}\left(  \operatorname{id}-z\square z\right)
^{-3/2}\left(  \operatorname*{d}\left(  z\square z\right)  \right)
z\right) \\
&  =&\frac{1}{2}m_{1}\left(  \left(  \operatorname{id}-z\square
z\right) ^{-2}z,\left(  \operatorname*{d}\left(  z\square z\right)
\right)  z\right) .
\end{eqnarray*}
We finally obtain , using this last result and (\ref{eq16}):

\begin{equation}
\Psi^{\ast}\left(  \overline{\partial}p_{1}\right)
=-\frac{\overline
{\partial}\mathcal{N}}{\mathcal{N}}+\frac{1}{2}m_{1}\left(  \left(
\operatorname{id}-z\square z\right)  ^{-2}z,\left(
\operatorname*{d}\left(
z\square z\right)  \right)  z\right)  . \label{eq18}%
\end{equation}

Along the same lines, one proves

\begin{equation}
\left(  \Psi^{-1}\right)  ^{\ast}\left(
\overline{\partial}p_{1}\right)
=\frac{\overline{\partial}\mathcal{N}_{*}}{\mathcal{N}_{*}}-\frac{1}{2}%
m_{1}\left(  \left(  \operatorname{id}+z\square z\right)
^{-2}z,\left(
\operatorname*{d}\left(  z\square z\right)  \right)  z\right)  . \label{eq19}%
\end{equation}

To prove (S) in Theorem \ref{mainteor}, it is now enough to check
that the forms
\begin{eqnarray*}
\beta(z)  &  =m_{1}\left(  \left(  \operatorname{id}-z\square
z\right) ^{-2}z,\left(  \operatorname*{d}\left(  z\square z\right)
\right)  z\right)
,\\
\beta_{*}(z)  &  =m_{1}\left(  \left(  \operatorname{id}+z\square
z\right) ^{-2}z,\left(  \operatorname*{d}\left(  z\square z\right)
\right)  z\right)
\end{eqnarray*}
are $\operatorname*{d}$-closed (or $\operatorname*{d}$-exact, as
$M$ and ${\cal M}$ are simply connected). We  verify it for $\beta
(z)$ in the following proposition (the proof for $\beta_*(z)$ is
similar).

\begin{prop}
Let $G$ be the analytic function defined on $\left]  -1,+1\right[
$ by
\[
G(t)=\frac{1}{t}\int_{0}^{t}\frac{u}{\left(  1-u\right)  ^{2}}%
\operatorname*{d}u
\]
and $\gamma:M\rightarrow\mathbb{R}$ the function defined by
\[
\gamma(z)=m_{1}\left(  G\left(  z\square z\right)  z,z\right)  .
\]
Then $\beta (z)=\operatorname{d}\gamma (z)$.
\end{prop}
\dimostr

By using Lemma \ref{lemmaa2}, in the appendix one has
\begin{equation*}
\operatorname*{d}\gamma(z)   =m_{1}\left(  G^{\prime}\left(
z\square z\right)  \left(  \operatorname*{d}\left(  z\square
z\right)  \right) z,z\right) +m_{1}\left(  G\left(  z\square
z\right)  \operatorname*{d}z,z\right) +m_{1}\left(  G\left(
z\square z\right)  z,\operatorname*{d}z\right)
\end{equation*}
Using the identity $G(t)+tG^{\prime}(t)=\frac{t}{\left(
1-t\right)  ^{2}},\label{eq22}$ we get
\begin{eqnarray*}
\operatorname*{d}\gamma(z) & =&m_{1}\left(  G^{\prime}\left(
z\square z\right)  \left(  \operatorname*{d}\left(  z\square
z\right)  \right)
z,z\right)  \\
&  -&m_{1}\left(  G^{\prime}\left(  z\square z\right)  \left(
z\square z\right)  \operatorname*{d}z,z\right)  -m_{1}\left(
G^{\prime}\left(
z\square z\right)  \left(  z\square z\right)  z,\operatorname*{d}z\right)  \\
&  +&m_{1}\left(  \left(  \operatorname{id}-z\square z\right)
^{-2}\left( z\square z\right)  \operatorname*{d}z,z\right)
+m_{1}\left(  \left( \operatorname{id}-z\square z\right)
^{-2}\left(  z\square z\right) z,\operatorname*{d}z\right)  .
\end{eqnarray*}
Now
\begin{eqnarray*}
 m_{1}\left(  G^{\prime}\left(  z\square z\right)  \left(  \operatorname*{d}%
\left(  z\square z\right)  \right)  z,z\right) & =&m_{1}\left(
G^{\prime}\left(  z\square z\right)  \left(
\operatorname*{d}z\square z\right)  z,z\right)  +m_{1}\left(
G^{\prime }\left(  z\square z\right)  \left(
z\square\operatorname*{d}z\right)
z,z\right)  \\
&  =&m_{1}\left(  G^{\prime}\left(  z\square z\right)  \left(
z\square z\right)  \operatorname*{d}z,z\right)  +m_{1}\left(
G^{\prime}\left( z\square z\right)  \left(
z\square\operatorname*{d}z\right)  z,z\right)
\end{eqnarray*}
and using the commutativity between $z\square z$ and $Q(z)$ and
the identity (4.55) in \cite{roos2}, p.495) one gets
\begin{eqnarray*}
m_{1}\left(  G^{\prime}\left(  z\square z\right)  \left(  z\square
\operatorname*{d}z\right)  z,z\right)   & =&m_{1}\left(
Q(z)\operatorname*{d}
z,G^{\prime}\left(  z\square z\right)  z\right)  \\
&  =&\overline{m_{1}\left(
\operatorname*{d}z,Q(z)G^{\prime}\left(  z\square
z\right)  z\right)  }\\
&  =&\overline{m_{1}\left(  \operatorname*{d}z,G^{\prime}\left(
z\square
z\right)  Q(z)z\right)  }\\
&=&m_{1}\left(  G^{\prime}\left(  z\square z\right)  \left(
z\square z\right)  z,\operatorname*{d}z\right).
\end{eqnarray*}
So we get%
\[
\operatorname*{d}\gamma(z)=m_{1}\left(  \left(
\operatorname{id}-z\square z\right)  ^{-2}\left(  z\square
z\right)  \operatorname*{d}z,z\right) +m_{1}\left(  \left(
\operatorname{id}-z\square z\right)  ^{-2}\left( z\square z\right)
z,\operatorname*{d}z\right)  .
\]
By the same argument as before, with $G^{\prime}$ replaced by
$F^{\prime}(t)=\frac{t}{(1-t)^2}$, we have
\begin{eqnarray*}
\beta (z)&=&m_{1}\left(  \left(  \operatorname{id}-z\square
z\right)  ^{-2}z,\left( \operatorname*{d}\left(  z\square z\right)
\right)  z\right)
 =m_{1}\left(  \left(  \operatorname{id}-z\square z\right)  ^{-2}\left(
z\square z\right)  \operatorname*{d}z,z\right) \\
&+&m_{1}\left(  \left( \operatorname{id}-z\square z\right)
^{-2}\left(  z\square z\right)
z,\operatorname*{d}z\right)=\operatorname{d}\gamma (z).
\end{eqnarray*}
\fdim

\section{Appendix: some technical results on HPJTS}

The following general result holds in Jordan triple systems (see
Lemma 2 in \cite{roosall}):

\begin{lemma}\label{lemmaa1}
Let ${\mathcal M}$ be an Hermitian positive Jordan triple system
with generic trace $m_{1}$ and generic norm $N$. Let
$\mathcal{N}(z)=N(z,z)$ and
$\mathcal{N}_{*}(z)=N(z,-z)$. Then%
\begin{eqnarray}
\frac{\overline{\partial}\mathcal{N}}{\mathcal{N}}  &
=-m_{1}\left(
z^{z},\operatorname*{d}z\right)  ,\label{eq12}\\
\frac{\overline{\partial}\mathcal{N}_{*}}{\mathcal{N}_{*}}  &
=m_{1}\left( z^{-z},\operatorname*{d}z\right),   \label{eq13}
\end{eqnarray}
where $z^z$ denotes the quasi inverse in the Jordan triple system
${\mathcal M}$.
\end{lemma}

\begin{lemma}\label{lemmaa2}
Let $f:\left]  -1,1\right[  \rightarrow\mathbb{R}$ and $F:\left]
-1,1\right[  \rightarrow\mathbb{R}$ be real-analytic functions.
Then
\begin{equation}
m_{1}\left(  f\left(  z\square z\right)  z,\left(
\operatorname*{d}F\left( z\square z\right)  \right)  z\right)
=m_{1}\left(  f\left(  z\square z\right)  z,F^{\prime}\left(
z\square z\right)  \operatorname*{d}\left(
z\square z\right)  z\right)  .\label{trace-gen}%
\end{equation}

\end{lemma}

\begin{proof}
It suffices to prove (\ref{trace-gen}) for $f=t^{p}$ and
$F=t^{k}$. For $k>0$, we have
\[
\left(  \operatorname*{d}\left(  \left(  z\square z\right)
^{k}\right)
\right)  z=\sum_{j=0}^{k-1}\left(  z\square z\right)  ^{k-1-j}%
\operatorname*{d}\left(  z\square z\right)  \left(  z\square
z\right)  ^{j}z.
\]
Recall that the \emph{odd powers} $z^{(2j+1)}$ in a Hermitian
Jordan triple
system are defined recursively by%
\begin{equation}
z^{(1)}=z,\quad z^{(2j+1)}=Q(z)z^{(2j-1)}\label{power}%
\end{equation}
and that they satisfy the identity
\begin{equation}
z^{(2j+1)}=\left(  z\square z\right)  ^{j}z.\label{eq29}%
\end{equation}
Using the commutativity between $z\square z$ and $Q(z)$ we then
have
\begin{eqnarray*}
\operatorname*{d}\left(  z\square z\right)  Q(z) &  =&\left(  \operatorname*{d}%
z\square z+z\square\operatorname*{d}z\right)  Q(z)=Q(z)\operatorname*{d}%
\left(  z\square z\right)  ,\\
\operatorname*{d}\left(  z\square z\right)  \left(  z\square
z\right)  ^{j}z
&  =&\operatorname*{d}\left(  z\square z\right)  Q(z)^{j}z=Q(z)^{j}%
\operatorname*{d}\left(  z\square z\right)  z,\\
\left(  \operatorname*{d}\left(  \left(  z\square z\right)
^{k}\right)
\right)  z &  =&\sum_{j=0}^{k-1}\left(  z\square z\right)  ^{k-1-j}%
Q(z)^{j}\operatorname*{d}\left(  z\square z\right)  z.
\end{eqnarray*}
As $z\square z$ is self-adjoint w.r. to $m_{1}$, we have
\begin{equation}
m_{1}\left(  z^{(2p+1)},\left(  \operatorname*{d}\left(  z\square
z\right) ^{k}\right)  z\right)  =\sum_{j=0}^{k-1}m_{1}\left(
\left(  z\square z\right)
^{p+k-1-j}z,Q(z)^{j}\operatorname*{d}\left(  z\square z\right)
z\right)  .\nonumber
\end{equation}
Using the identity (4.55) in \cite{roos2}, p.495, we obtain
(denoting by
$\tau$ the conjugation of complex numbers)%
\begin{eqnarray*}
 m_{1}\left(  z^{(2p+1)},\left(  \operatorname*{d}\left(  z\square z\right)
^{k}\right)  z\right) &=&\sum_{j=0}^{k-1}\tau^{j}m_{1}\left(
Q(z)^{j}\left(  z\square z\right)
^{p+k-1-j}z,\operatorname*{d}\left(  z\square z\right)  z\right)  \\
&  =&\sum_{j=0}^{k-1}\tau^{j}m_{1}\left(  \left(  z\square
z\right)
^{p+k-1}z,\left(  z\square z\right)  \operatorname*{d}z+Q(z)\operatorname*{d}%
z\right)  .
\end{eqnarray*}
But
\begin{eqnarray*}
m_{1}\left(  \left(  z\square z\right)
^{p+k-1}z,\operatorname*{d}\left( z\square z\right)  z\right)
 &=&m_{1}\left(  \left(  z\square z\right)  ^{p+k-1}z,\left(  z\square
z\right)  \operatorname*{d}z+Q(z)\operatorname*{d}z\right)  \\
 &=&m_{1}\left(  z^{(2p+2k+1)},\operatorname*{d}z\right)  +\tau m_{1}\left(
z^{(2p+2k+1)},\operatorname*{d}z\right)
\end{eqnarray*}
is real, so that
\begin{eqnarray*}
m_{1}\left(  z^{(2p+1)},\left(  \operatorname*{d}\left(  z\square
z\right) ^{k}\right)  z\right)   &  =&km_{1}\left(  \left(
z\square z\right)
^{p+k-1}z,\operatorname*{d}\left(  z\square z\right)  z\right)  \\
&  =&m_{1}\left(  \left(  z\square z\right)  ^{p}z,k\left(
z\square z\right) ^{k-1}\operatorname*{d}\left(  z\square z\right)
z\right)  ,
\end{eqnarray*}
which is precisely (\ref{trace-gen}) for $f=t^{p}$, $F=t^{k}$.
\end{proof}

\small{}


\begin{thebibliography}{99}

\bibitem{Ber} W. Bertram,
{\em The geometry of Jordan and Lie structures}, Lecture Notes in
Mathematics 1754, Springer-Verlag (2000).


\bibitem{cr1} E. Ciriza,
{\em The local structure of a Liouville vector field}, Amer. J.
Math. 115 (1993), 735-747.

\bibitem{cr2} E. Ciriza,
{\em On special submanifolds in symplectic geometry}, Diff. Geom.
Appl. 3 (1993), 91-99.

\bibitem{cr3} E. Ciriza,
{\em Symplectomorphic codimension $1$ totally geodesic
submanifolds}, Diff. Geom. Appl. 5 (1995), 99-104.



\bibitem{alcu} F. Cuccu and A. Loi,
{\em Global symplectic coordinates on complex domains}, to appear
in Journal of Geometry and Physics.

\bibitem{c1nash} G. D'Ambra and A. Loi
{\em A symplectic version of Nash $C^1$-isometric embedding
theorem}, Diff. Geom. Appl. 16, no. 2 (2002), 167-179.

\bibitem{connmet} G. D'Ambra and A. Loi,
{\em Inducing connections on SU(2)-bundles}, JP Journal of
Geometry and Topology 3 (2003), 65-88.


\bibitem{gr} M. Gromov,
{\em Partial Differential Relations}, Springer-Verlag (1986).

\bibitem{JaPa} N. Jacobson and L. Page,
{\em  On Jordan algebras with two generators}, J. Math. Mech. 6
(1957), 895--906.



\bibitem{kn}S. Kobayashi and K. Nomizu,
{\em Foundations of Differential Geometry vol. II}, John Wiley and
Sons Inc. (1967).

\bibitem{ko1} M. Koecher,
{\em The Minnesota Notes on Jordan Algebras and Their
Applications}, Lecture Notes in Mathematics 1710, Springer-Verlag
(1999).


\bibitem{ko2} M. Koecher,
{\em An elementary approach to Bounded Symmetric Domains,} Rice
University (1969).


\bibitem{diastherm} A. Loi,
{\em Calabi's diastasis function for Hermitian symmetric spaces},
to appear in Diff. Geom. Appl. available in:
\mbox{http://loi.sc.unica.it/articoli.html - no. 18}.


\bibitem{loos} O. Loos,
{\em Bounded Symmetric Domains and Jordan pairs,} Lecture Notes,
Irvine (1977).

\bibitem{mc} D. McDuff,
{\em The symplectic structure of \K\ manifolds of non-positive
curvature}, J. Diff. Geometry 28  (1988), 467-475.

\bibitem{mok} N. Mok,
{\em Metric Rigidity Theorems on Hermitian Locally Symmetric
Spaces}, Series in Pure Mathematics -Volume 6. World Scientific
(1989).


\bibitem{roos} G. Roos,
{\em Volume of bounded symmetric domains and compactification of
Jordan triple systems}, Lie groups and Lie algebras, 249--259,
Math. Appl., 433, Kluwer Acad. Publ., Dordrecht, 1998.

\bibitem {roos2}G. Roos, Jordan triple systems, pp. 425-534, in
\textit{J.~Faraut, S.~Kaneyuki, A.~Kor\'{a}nyi, Q.k.~Lu, G.~Roos,
Analysis and Geometry on Complex Homogeneous Domains}, Progress in
Mathematics, vol.\textbf{185}, Birkh\"{a}user, Boston, 2000.


\bibitem{roosall} G. Roos et all.,
{\em Kaehler--Einstein metric for some Hartogs domains over
bounded symmetric domains}, ArXiv math.CV/0501223, 2005.

\bibitem{sa} I. Satake,
{\em Algebraic structures of symmetric domains}, Publications of
the Mathematical Society of Japan 14, Kano Memorial Lectures 4,
Iwanami Shoten Pub. and Princeton University Press (1980).


\bibitem{ta} M. Takeuchi,
{\em Homogeneous \K\ Manifolds in Complex Projective Space}, Japan
J. Math. vol. 4 (1978), 171-219.



\bibitem{wolf} J. A.  Wolf,
{\em Fine structure of Hermitian symmetric spaces}, Pure and App.
Math. 8 (1972), 217-357.





\end{thebibliography}
\end{document}